\documentclass[12pt]{amsart}
\usepackage{amsfonts,amssymb,amsmath,amsthm,amscd,graphicx,mathtools, float}
\DeclarePairedDelimiter\floor{\lfloor}{\rfloor}
\newtheorem{theorem}{Theorem}[section]
\newtheorem{lemma}[theorem]{Lemma}
\newtheorem{proposition}[theorem]{Proposition}

\newtheorem{corollary}[theorem]{Corollary}

\theoremstyle{definition}

\theoremstyle{remark}
\newtheorem{remark}[theorem]{Remark}

\numberwithin{equation}{section}

\newcommand{\pZ}{\mathbb{Z}_{\geq 0}}

\newcommand{\lcr}{\raisebox{-5pt}{\mbox{}\hspace{1pt}
                 \includegraphics{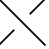}\hspace{1pt}\mbox{}}}
\newcommand{\ift}{\raisebox{-5pt}{\mbox{}\hspace{1pt}
                 \includegraphics{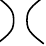}\hspace{1pt}\mbox{}}}
\newcommand{\zer}{\raisebox{-5pt}{\mbox{}\hspace{1pt}
                 \includegraphics{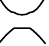}\hspace{1pt}\mbox{}}}
\newcommand{\rcr}{\raisebox{-5pt}{\mbox{}\hspace{1pt}
                 \includegraphics{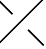}\hspace{1pt}\mbox{}}}

\title{The Localized Skein Algebra is Frobenius}

\author{Nel Abdiel}

\author{Charles Frohman}

\address{Department of Mathematics, The University of Iowa, Iowa City, IA
52242, USA}

\begin{document}

\begin{abstract} When $A$ in the Kauffman bracket skein relation is a primitive $2N$th root of unity, where $N\geq 3$ is odd, the Kauffman bracket skein algebra $K_N(F)$ of a finite type surface $F$ is a ring extension of the $SL_2\mathbb{C}$-characters $\chi(F)$  of the fundamental group of $F$. We localize by inverting the nonzero characters to get an algebra $S^{-1}K_N(F)$ over the function field of the character variety.  We prove that if $F$ is noncompact, the algebra $S^{-1}K_N(F)$ is a symmetric Frobenius algebra. Along the way we prove $K(F)$  is finitely generated, $K_N(F)$ is a finite rank module over $\chi(F)$, and the simple closed curves that make up any simple diagram on $F$ generate a finite field extension of $S^{-1}\chi(F)$ inside $S^{-1}K_N(F)$. \end{abstract}

\maketitle

\section{Introduction}

 This paper is a step in a program to build a $4$-dimensional extended field theory that assigns invariants to manifolds equipped with a homomorphism of their fundamental group into $SL_2\mathbb{C}$. Frobenius algebras are central to the construction of field theories. In this paper we show that the Kauffman bracket skein algebra of a compact surface with nonempty boundary can be localized to give a symmetric Frobenius algebra over the function field of the $SL_2\mathbb{C}$-character variety of the fundamental group of the surface.
 
 A surface $F$ is of finite type if there is a closed oriented surface $\hat{F}$ and a finite set of points $\{p_j\}\in \hat{F}$ so that $F=\hat{F}-\{p_i\}$. In this paper all surfaces are either compact  oriented (possibly with boundary) or of finite type.  If $F$ is a compact connected oriented surface, a punctured disk can be glued into each boundary component to obtain a finite type surface. There is a one-to-one correspondence between disjoint families of simple closed curves in the two surfaces, so we can prove our theorems by only working with finite type surfaces.

We  prove that $K(F)$ the Kauffman bracket skein algebra over $\mathbb{Z}[A,A^{-1}]$ is finitely generated as an algebra by a finite family of simple closed curves $S_i$. In fact, \[ \{S_{\sigma(1)}^{k_1}*S_{\sigma(2)}^{k_2}*\ldots* S_{\sigma(n)}^{k_n}\},\] where $k_i\in \mathbb{Z}_{\geq 0}$, spans $K(F)$ for any permutation $\sigma:\{1,\ldots,n\}\rightarrow \{1,\ldots,n\}$. This is an extension of a theorem of Bullock \cite{Bu}.
Key to the proof is the description of a graded algebra derived from a filtration of $K(F)$ when $F$ is not closed. The filtration comes via geometric intersection numbers with the edges of an ideal triangulation.  Geometric intersection numbers have been used as a tool to filter the skein algebra to good effect in \cite{Le,M,Mu}. In \cite{Mu}, Muller studies a larger algebra where there are many more filtrations based on geometric intersections that can be applied.

Denote the result of specializing $A$ at a primitive $2N$th root of unity, $N\geq 3$ and odd, by $K_N(F)$.
We use the result above to prove that $K_N(F)$ is a finitely generated module over $\chi(F)$, the coordinate ring of the $SL_2\mathbb{C}$-character variety of $\pi_1(F)$.  Localizing at $S=\chi(F)-\{0\}$, we get $S^{-1}K_N(F)$ is a finite dimensional algebra over $S^{-1}\chi(F)$, the function field of the character variety of $\pi_1(F)$.

Each $\alpha \in S^{-1}K_N(F)$ induces an $S^{-1}\chi(F)$-linear endomorphism \[L_{\alpha}:S^{-1}K_N(F)\rightarrow S^{-1}K_N(F),\] by left multiplication. The {\bf normalized trace} of $\alpha$, $Tr(\alpha)$ is the trace  of $L_{\alpha}$ as a linear endomorphism divided by the dimension of the vector space $S^{-1}K_N(F)$ over the field $S^{-1}\chi(F)$.  The normalized trace has the properties,
\begin{itemize}
\item $Tr(1)=1$,
\item $Tr(\alpha*\beta)=Tr(\beta*\alpha)$
\item $Tr$ is $S^{-1}\chi(F)$-linear.

\end{itemize}  Hence if $Tr$ is nondegenerate then $S^{-1}K_N(F)$ equipped with the normalized trace is a symmetric Frobenius algebra over the function field of the character variety of the fundamental group of the surface $F$.

Along the way we learn to compute $Tr:S^{-1}K_N(F)\rightarrow S^{-1}\chi(F)$
with respect to a special basis. A primitive diagram on $F$ is a system of disjoint simple closed curves $S_i$ so that no $S_i$ bounds a disk and no two curves in the system cobound an annulus. The skein $\prod_iT_{k_i}(S_i)$ is the product over all $i$ of the result of threading $S_i$ with the $k_i$th Tchebychev polynomial of the first kind.   These span $S^{-1}K_N(F)$ over $S^{-1}\chi(F)$.  If
\[ \alpha=\sum_{k_1,\ldots,k_n}\alpha_{k_1,\ldots,k_n}\prod_iT_{k_i}(S_i)\]
where the $\alpha_{k_1,\ldots,k_n}\in S^{-1}\chi(F)$ then $Tr(\alpha)$ is the sum of all the terms where $N|k_i$ for all $i$. The derivation of the formula for the trace depends on the following surprising fact.  Let $\cup_iS_i$ be a simple diagram, made up of the simple closed curves $S_i$. The extension of $S^{-1}\chi(F)$ obtained by adjoining the $S_i$, $S^{-1}\chi(F)[S_1,\ldots,S_n]$ is a field.  This extends a result of Muller \cite{Mu}, that says that simple closed curves are not zero divisors.

Since the value of the formula for the trace doesn't have any denominators that didn't appear in the input, the trace is actually defined as a $\chi(F)$-linear map
\[ Tr:K_N(F)\rightarrow \chi(F).\]

Next, the formula for the trace is used to prove that there are no non-trivial principal ideals in the kernel of $Tr:K_N(F)\rightarrow \chi(F)$, completing the proof that $S^{-1}K_N(F)$ is a symmetric Frobenius algebra.  Essential to the proof is the fact that given  a primitive diagram $\cup S_i$ ,  the skeins $\prod_iT_{k_i}(S_i)$ with $k_i\geq 0$ generate a field extension of $S^{-1}\chi(F)$ in $S^{-1}K_N(F)$ . 

The authors would like to thank Pat Gilmer and Thang L\^{e} for helpful input.

\section{Preliminaries}
\subsection{Kauffman Bracket Skein Module}
Let $M$ be an orientable $3$-manifold.  A {\bf framed link} in $M$ is an
embedding of a disjoint union of annuli into $M$.   Throughout this paper $M=F\times [0,1]$ for an orientable surface $F$. Diagrammatically we  depict framed links by showing the core of the annuli  lying parallel to $F$.  Two framed links in $M$ are equivalent if they are isotopic. Let $\mathcal{L}$ denote the set of
equivalence classes of framed links in $M$, including the empty
link. By $\mathbb{Z}[A,A^{-1}]$ we mean Laurent polynomials with integral coefficients in the formal variable $A$.
Consider the free module over  $\mathbb{Z}[A,A^{-1}]$,
\[\mathbb{Z}[A,A^{-1}]\mathcal{L}\] with basis $\mathcal{L}$.  Let $S$ be the submodule spanned by the Kauffman bracket skein relations, 

\[\lcr-A\zer-A^{-1}\ift \]
and
\begin{equation*}\bigcirc \cup L+(A^2+A^{-2})L.\end{equation*} 
 The framed links in 
each expression are
identical outside the balls pictured in the diagrams, and when both arcs in a diagram lie in the same component of the link, the same side of the annulus is up.  
 The Kauffman bracket
skein module $K(M)$ is the quotient
\[ \mathbb{Z}[A,A^{-1}] \mathcal{L} / S(M). \]
A {\bf skein} is an element of $K(M)$.
Let $F$ be a compact orientable surface and let $I=[0,1]$. There is an
algebra structure on $K(F\times I)$ that comes from laying one framed link
over the other.  The resulting algebra is denoted
$K(F)$ to emphasize
that it comes from the particular structure as a cylinder over $F$. Denote the stacking product with a $*$, so $\alpha*\beta$ means $\alpha$ stacked over $\beta$. If it is known the two skeins commute  the $*$ will be omitted.

A {\bf simple diagram} $D$ on the surface $F$ is a system of disjoint simple closed curves so that none of the curves bounds a disk.    A simple diagram $D$ is {\bf primitive} if no two curves in the diagram cobound an annulus.  A simple diagram can be made into a framed link by choosing a system of disjoint annuli in $F$ so that each annulus has a single curve in the diagram as its core. This is sometimes called the {\bf blackboard framing}.  The set of isotopy classes of blackboard framed simple diagrams form a basis for $K(F)$, \cite{BFK,HP,SW}.

\subsection{Specializing $A$.}  If $R$ is a commutative  ring , and $\zeta \in R$ is a unit, then $R$
 a $\mathbb{Z}[A,A^{-1}]$-module, where the action 
\[ \mathbb{Z}[A,A^{-1}]\otimes  R \rightarrow R\]
is given by letting $p\in \mathbb{Z}[A,A^{-1}]$ act by multiplication by the result of evaluating $p$ at $\zeta$.  The skein module specialized at $\zeta \in R$ is,
\[ K_R(M)=K(M)\otimes_{\mathbb{Z}[A,A^{-1}]}R.\]
You can think of the specialization as setting $A$ equal to $\zeta$ in the Kauffman bracket skein relations.

This is much too general a setting to get nice structure theorems for $K_R(M)$, so
we restrict our attention to when the ring $R$ is an integral domain. To emphasize that we are working with an integral domain we denote the ring $\mathfrak{D}$.  Since
$\mathbb{Z}[A,A^{-1}]$ is an integral domain and $A$ is a unit we can recover $K(M)$ by specialization. For that reason the theorems in this paper are all stated in terms of $K_{\mathfrak{D}}(M)$, the skein module specialized at a unit $\zeta$ in an integral domain $\mathfrak{D}$.

We are most interested in the case when $\zeta$ is a primitive $2N$th root of unity, where $N\in \pZ$ is odd. The integral domain is $\mathbb{Z}[\frac{1}{2},\zeta]$ which we think of as embedded in $\mathbb{C}$, so that $\zeta=e^{k\pi{\bf i}/N}$, where $k$ is an odd counting number that is relatively prime to $N$.  In this case, the specialized module is denoted $K_N(M)$ .  This is a little ambiguous because there are in general several primitive $2N$th roots of unity, but the theorems in this paper do not depend on the choice of which one. We need $2$ to be a unit so that a collection of skeins that are adapted to the computation of the trace will be a basis.

\subsection{Threading}

The Tchebychev Polynomials of the first type  $T_k$ are defined recursively by
\begin{itemize}
\item $T_0(x) =2$, \item $T_1(x)= x$ and \item $T_{n+1}(x) = T_1(x) \cdot T_n(x) - T_{n-1}(x)$. \end{itemize}

They satisfy some nice properties.
\begin{proposition}   \label{mult} 
For $m,n>0$, $T_m(T_n(x))=T_{mn}(x)$. Furthermore, for all $m,n\geq 0$,  $T_m(x)T_n(x)=T_{m+n}(x)+T_{|m-n|}(x)$.\end{proposition}
For a proof see \cite{AF}.

Let $\Sigma_{0,2}=S^1\times [0,1]$. It is easy to see that $K_{\mathfrak{D}}(\Sigma_{0,2})$ is isomorphic to $\mathfrak{D}[x]$
where $x$ is the framed link coming from the blackboard framing of the core of the annulus.  Hence $1,x,x^2,\ldots x^n,\ldots$ is a basis for $K_{\mathfrak{D}}(\Sigma_{0,2})$.  Since $T_0(x)=2$, in order to used the $T_k(x)$ as a basis for $\mathfrak{D}$, $2$ must be a unit  in $\mathfrak{D}$.  If $2\in \mathfrak{D}$ is a unit then $\{T_k(x)\}$, $k\in \pZ$ is a basis for $K_{\mathfrak{D}}(\Sigma_{0,2})$.

If the components of the primitive diagram on a finite type surface $F$ are the  simple closed curves $S_i$  and $k_i\in \mathbb{Z}_{\geq 0}$
has been chosen for each component, the result of  threading the diagram with the $T_{k_i}$ is $\prod_iT_{k_i}(S_i)$.  Since the $S_i$ are disjoint from one another they commute so order doesn't matter in the product.
For any compact or finite type surface $F$, the primitive diagrams on $F$ up to isotopy, with their components threaded with $\{T_{k_i}\}$ form a basis for $K_{\mathfrak{D}}(F)$ so long as $2\in \mathfrak{D}$ is a unit. This basis is becoming more commonly used in the study of skein algebras \cite{FG,Th,Le2}.

The following theorem of Bonahon and Wong is the starting point for this investigation.

\begin{theorem}
(Bonahon-Wong \cite{BW,Le}). If  $M$ is a compact oriented three-manifold and we specialize at $\zeta$ a $2N$th root of unity with $N\geq 3$ odd, there is a $\mathbb{Z}[\frac{1}{2},\zeta]$-linear map
\begin{center}
$\tau: K_{1}(M) \rightarrow K_{N}(M)$
\end{center}
given by threading framed links with $T_N$. Any framed link in the image of $\tau$ is central in the sense that if $L' \cup K$ differs from $L \cup K$ by a crossing change of $L$ and $L'$ with $K$, then $T_N(L) \cup K = T_N(L') \cup K$. In the case that $M = F \times [0,1]$, the map
\begin{center}
$\tau: K_{1}(F) \rightarrow K_N(F)$
\end{center}
is an injective homomorphism of algebras so that the image of $\tau$ lies in the center of $K_N(F)$.
\end{theorem}

The skein module $K_1(M)$ is a ring under disjoint union. At $A=-1$, the Kauffman bracket skein relation 
\[\lcr+\zer+\ift \]
 can be rotated $90$ degrees and then subtracted from itself  to yield,
 \[ \lcr-\rcr.\]  This means that in $K_1(M)$ changing crossings does not change the skein. To take the product of two equivalence classes of framed links, choose representatives that are disjoint from one another and take their union. The product
is independent of the representatives chosen, since the results differ by isotopy and changing crossings.
The product can be extended distributively to give a product on $K_1(M)$. Let $\sqrt{0}$ denote the nilradical of $K_1(M)$.  It is a theorem of Bullock, \cite{Bu}, proved independently in \cite{PS}, that for any oriented compact $3$-manifold $K_1(M)/\sqrt{0}$ is canonically isomorphic to the coordinate ring of the $SL_2\mathbb{C}$-character variety of the fundamental group of $M$.  In the case that $M=F\times [0,1]$ the disjoint union product coincides with the stacking product, as stacking is one way to perturb the components of the two links so that they are disjoint. It is a theorem of Przytycki and Sikora \cite{PS} that  $\sqrt{0}=\{0\}$ when the underlying three-manifold is $F\times[0,1]$.  Therefore, $K_1(F)$ is the coordinate ring of the $SL_2\mathbb{C}$-character variety of the fundamental group of $F$. To alleviate notational ambiguity, and to emphasize the relationship with character varieties,  the image of the threading map is denoted by $\chi(M)$.

For any oriented finite type surface $F$, $\chi(F)$ has basis  the isotopy classes of primitive diagrams  threaded with $T_{k_iN}$ for all choices $k_i\in \mathbb{Z}_{\geq 0}$.

\subsection{Specializing at a place.}  A place of $\chi(F)$ is a homomorphism $\phi:\chi(F)\rightarrow \mathbb{C}$.  The places correspond to evaluation at a point on the character variety. 
A place defines a module structure $\chi(F)\otimes \mathbb{C}\rightarrow \mathbb{C}$
by letting $s\in \chi(F)$ act as multiplication by $\phi(s)$ on $\mathbb{C}$.  We define the {\bf specialization} of $K_N(F)$ at $\phi$ to be,
\[ K_N(F)_{\phi}=K_N(F)\otimes_{\chi(F)} \mathbb{C}.\]
The specialization at a place is an algebra over the complex numbers. 

\subsection{Localization}

Let $R\rightarrow J$ be a central ring extension, where is $R$ is an integral domain, $J$ is an associative ring with unit
and the inclusion of $R$ into $J$ is a ring homomorphism.  Since $R$ has no zero divisors, $S=R-\{0\}$ is multiplicatively closed.   Start with the set of ordered pairs $J\times S$, and  place an equivalence relation on $J\times S$ by  saying $(a,s)$  is equivalent to $(b,t)$ if $at=bs$.  Denote the equivalence class of $(a,s)$ under this relation by $[a,s]$.  The set of equivalence classes is denoted $S^{-1}J$, and called the localization of $J$ with respect to  $S$. Denote the set of equivalence classes $[a,s]$ where $a\in R$ by $S^{-1}R$. Define multiplication of equivalence classes by $[a,s][b,t]=[ab,st]$ and addition by $[a,s]+[b,t]=[at+bs,st]$.  Under these operations $S^{-1}R$ is a field, and $S^{-1}J$ is an algebra over that field.

In this paper, $J$ is  a subalgebra of $K_N(F)$ and $R$ is $\chi(F)$.  This means that $S^{-1}R$ is the function field of the character variety of $\pi_1(F)$.

\subsection{Trace and Extension of Scalars}

If $L\in End_F(V)$, we use $tr(L)$ to denote the unormalized trace of $L$.    The linear map $L$ can be represented with respect to a basis $\{v_j\}$, by a matrix $(l^j_i)$. The trace of $L$ is given by \[tr(L)=\sum_il^i_i.\] If $W$ is also a finite dimensional vector space over $F$ and $M:W\rightarrow W$ is an $F$-linear map, then
\[ tr(L\otimes_F M)=tr(L)tr(M) \ \mathrm{and} \ tr(L\oplus M)=tr(L)+tr(M).\]

Suppose that $F\leq K$ is a field extension and $V$ is a vector space of dimension $n$ over $F$, then 
\[ V\otimes_F K\]
is a vector space of dimension $n$ over $K$. In fact if $\{v_j\}$ is a basis for $V$ then $\{v_j\otimes 1\}$ is a basis for $V\otimes_F K$ over $K$.

Under extension of scalars, $L:V\rightarrow  V$ gets sent to $L\otimes_F 1_K$. The matrix of  $L\otimes_F 1_K$ with respect to the basis $\{v_j\otimes 1\}$ is the same as the matrix of $L$ with respect to $\{v_j\}$, so
\[tr(L\otimes_F 1_K)=tr(L),\] where the trace on the left is taken as a $K$-linear map,
and the trace on the right is taken as an $F$-linear map, and we are using $F\leq K$ to make the identification.

\subsection{Geometric Intersection Numbers}

Suppose that $X$ and $Z$ are properly embedded $1$-manifolds in the finite type surface $F$ where  $X$ is compact.   We say that $X'$ is a transverse representative of $X$, if $X'$ is ambiently isotopic to $X$ via a compactly supported isotopy, and $X'\pitchfork Z$.   Define the {\bf geometric intersection number} of $X$ and $Z$, denoted $i(X,Z)$ to be the minimum cardinality of $X'\cap Z$ over all transverse representatives of $X$. We could have instead worked with $Z$ up to compactly supported ambient isotopy and taken the minimum over all $Z'$ isotopic to $Z$ and transverse to $X$ and gotten the same number, so $i(X,Z)=i(Z,X)$.  

It is a theorem that a transverse representative of $X$ realizes the geometric intersection number $i(X,Z)$ if and only if there are no {\bf bigons}. A bigon is a disk $D$ embedded in $F$ so that the boundary of $D$ consists of the union of two arcs $a\subset X$ and $b\subset Z$, \cite{FLP}. If there is a bigon, there is always an innermost bigon, whose interior is disjoint from $X\cup Z$.

\section{$K_{\mathfrak{D}}(F)$ is  finitely generated.}

\subsection{ Parametrizing the simple diagrams.}

An ideal triangle is a triangle with its vertices removed.
An ideal triangulation of a finite type surface $F$ consists of finitely many  ideal triangles $\Delta_i$ with their edges identified pairwise, along with a homeomorphism from the resulting quotient space to $F$. Alternatively an ideal triangulation is defined by a family $C$ of properly embedded lines that cuts $F$ into finitely many  ideal triangles. More precisely, place a complete Riemannian metric on $F$. Define a metric on the components $D$ of $F-C$ by
\[ d(p,q)=\inf{\{length(\alpha)|\alpha:[0,1]\rightarrow D \ \mathrm{is \ smooth}, \alpha(0)=p,\alpha(1)=q}\}.\] 
Let $\Delta$ be the metric space completion of $D$. If $\Delta$ is homeomorphic to a disk with three points removed from its boundary, then $D$ {\bf completes to an ideal triangle}.
If all the components of $F-C$ complete to ideal triangles then $C$ {\bf cuts $F$ into ideal triangles.}

If $\Delta$ is an ideal triangle in an ideal triangulation then  $\partial \Delta=\{a,b,c\}$ where $a$, $b$ and $c$ are homeomorphic to $\mathbb{R}$. The lines $a$, $b$, and $c$ are the {\bf sides} of $\Delta$. 
There is a map of $\Delta$ to the closure of a component $D$ of the complement of $C$ into $F$.  If this map is an embedding, then $\Delta$ is an {\bf embedded ideal triangle}.  It could be that two sides $c_1, c_2$ of the ideal triangle $\Delta$ get mapped to the same line $c$, in this case $\Delta$ is a {\bf folded ideal triangle}.   Figure \ref{folded} is a picture of a folded ideal triangle.  There are two punctures in the picture, and the mapping is $2-1$ along the vertical line joining them. The edge that is covered twice by the mapping has {\bf multiplicity} $2$.

\begin{figure}[h] \label{folded} \includegraphics{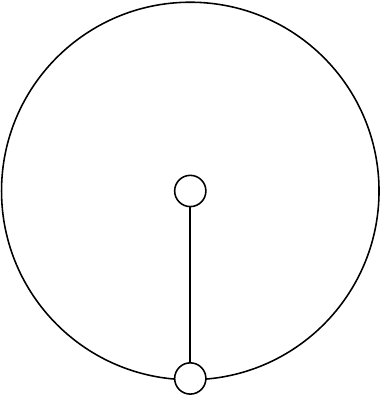}\caption{ A folded triangle}\end{figure}

 The surface $F$ needs to have at least one puncture,  and  negative Euler characteristic or it doesn't admit an ideal triangulation.
If the Euler characteristic of the surface $F$ is $-e(F)$ then any ideal triangulation of $F$ consists of  $2e(F)$ ideal triangles. The cardinality of a set of lines $C$ defining an ideal triangulation is $3e(F)$.
 
Let $C$ denote a disjoint family of properly embedded lines that defines an  ideal triangulation of $F$, and suppose the triangles are the set $\{\Delta_j\}$. 
An {\bf admissible coloring},  $f:C\rightarrow \pZ$ is an assignment of a nonnegative integer $f(c)$ to each $c\in C$ so that the following conditions hold;
\begin{itemize}
\item If $\{a,b,c\}$ form the boundary of an embedded ideal triangle  $\Delta_j$ then $f(a)+f(b)+f(c)$ is even and the triple $\{f(a),f(b),f(c)\}$ satisfies the triangle inequality,
\[ f(a)\leq f(b)+f(c), \ f(b)\leq f(a)+f(c), \ \mathrm{and} \  f(c)\leq f(a)+f(b).\]
\item If $\{a,b\}$  are the image of the boundary of  a folded ideal triangle $\Delta_j$  where $b$ has multiplicity $2$ we require that $f(a)+2f(b)$ be even and $f(a)\leq 2f(b)$.
\end{itemize}
If $S\subset F$ is a simple diagram then $f_S:C\rightarrow \mathbb{Z}_{\geq 0}$ given by $f_S(c)=i(S,c)$ is an admissible coloring.

\begin{proposition} The admissible colorings $f:C\rightarrow \mathbb{N}$ are in one-to-one correspondence with isotopy classes of simple diagrams on $F$.\end{proposition}

\qed 

The admissible colorings form a {\bf pointed integral polyhedral cone}. If $C$ is the set of edges of the ideal triangulation then there is a map,
\[ \{f:C\rightarrow \pZ|admissible\} \rightarrow \mathbb{Z}^C\]
that sends each $f$ to its tuple of values.  The image is defined by linear equations and inequalities, so it is polyhedral. The image is closed under addition, so it is a cone, and the tuple of all zeroes is an admissible coloring so it is pointed.

It is a classical result \cite{G} that any pointed integral polyhedral cone admits an {\bf integral basis}.  That is, there are finitely many admissible colorings $f_{S_i}$ so that every admissible coloring is a nonnegative integral linear combination of the $f_{S_i}$, and the set $\{f_{S_i}\}$ has minimal cardinality with respect to this condition.  The integral basis is unique.
If $P$ is a pointed integral polyhedral cone, $p\in P$ is {\bf indivisible} if whenever $s,t\in P$ and $s+t=p$ then $s=0$ or $p=0$. The set of indivisible elements of $P$ is the integral basis \cite{V}.  In the case of the cone of admissible colorings, the diagrams corresponding to indivisible colorings are simple closed curves.

 Forgetting positivity, and the triangle inequality, the admissible colorings generate a free module over $\mathbb{Z}$. It makes sense to ask whether a collection $f_{S_i}:C\rightarrow \pZ$ are linearly independent. Oddly, the integral basis need not be linearly independent.

\subsection{An Example}
Decompose the punctured torus $\Sigma_{1,1}$ into two ideal triangles. This requires three edges, that form the boundary of both triangles. In the diagram below  we identify the left and right hand sides of the rectangle, and the top and bottom of the rectangle  with the vertices deleted to obtain a once punctured torus. The lines defining the triangulation come from the sides of the rectangle and the diagonal shown. 

\begin{figure}[h] \includegraphics{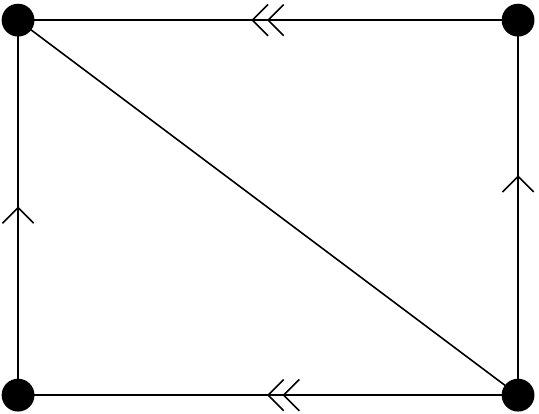}\caption{ An ideal triangulation of $\Sigma_{1,1}$}\end{figure}

The admissible colorings can be seen as triples of counting numbers $(m,n,p)$ whose sum is even and satisfy the triangle inequality. The nonzero indecomposable admissible colorings are $(1,1,0),(1,0,1)$ and $(0,1,1)$. This set is an integral basis.  Notice that if $(a,b,c)$ is an admissible coloring and one of the triangle inequalities is strict, say $a<b+c$ we can subtract the corresponding indecomposable $(0,1,1)$ to get a triple $(a,b-1,c-1)$  that still satisfies the triangle inequality and  the sum of the colors $a+b+c-2<a+b+c$. If all three triangle inequalities are equalities, $a=b+c$, $b=a+c$ and $c=a+b$ then $(a,b,c)=(0,0,0)$.
The three curves corresponding to $(1,1,0)$, $(1,0,1)$ and $(0,1,1)$  are the generators that Bullock and Przytycki \cite{BP}  got for $K(\Sigma_{1,1})$. There are infinitely many ideal triangulations of $\Sigma_{1,1}$ but Euler characteristic forces them all to be two triangles that share all their edges. The argument above goes through, even though the curves on the torus will be different. Since the  integral basis is unique, any set of skeins that generates $K_{\mathfrak{D}}(\Sigma_{1,1})$ must have at least three elements.

\subsection{The algebra $K_{\mathfrak{D}}(F)$ is finitely generated over $\mathfrak{D}$.}

If $f_S$ and $f_{S'}$ are admissible colorings, choose simple diagrams $S$ and $S'$ that realize the colorings as the cardinality of their intersections with the $c_i \in C$ and so that $S$ and $S'$ realize their geometric intersection number, $S\cap S'$ is disjoint from all $c_i$.  Up to isotopy there is a unique simple diagram whose associated coloring is $f_S+f_{S'}$,
called the {\bf geometric sum} of $S$ and $S'$.  Since addition of admissible colorings is associative, so is the geometric sum.

More is true. Define the {\bf weight} of a simple diagram $S$ to be \[weight(S)=\sum_{c\in C}i(S,c)=\sum_{c\in C}f_S(c).\] Suppose that $S$ and $S'$ transversely represent $i(S,S')$. Furthermore assume that $S\cap S' \cap C=\emptyset$.  If there are $n$ points of intersection in $S\cap S'$ there are $2^n$ ways of smoothing all the crossings of $S$ and $S'$ to get a system of simple closed curves.  We call a system of simple closed curves obtained by smoothing all crossings $s$ a {\bf state}. A state might not be a simple diagram as it may contain some trivial simple closed curves. There is a process for writing the product $S*S'$ as a linear combination of simple diagrams.  First expand the product as a sum of states using the Kauffman bracket skein relation for crossings, then delete the trivial components of each state, and for each trivial component deleted from a state multiply the coefficient of the state by $-\zeta^2-\zeta^{-2}$. Order the crossings of $S*S'$. Based on the ordering there is a rooted tree, where the root is the diagram $S*S'$, the vertices are partial smoothings (resolvents) of the diagram, and the directed edges correspond to smoothing the crossings in order.  The states are the leaves of this tree. If the shortest path from the root to a state $s$ passes through a resolvent $r$ we say that $s$ is a {\bf descendent} of $r$.

\begin{theorem} \label{weighty} Let $S$ and $S'$ be simple diagrams associated to admissible colorings $f_S,f_{S'}:C\rightarrow \mathbb{Z}_{\geq 0}$.  Assume the product $S*S'\in K_{\mathfrak{D}}(F)$ has been written  as $\sum_D\alpha_D D$ where the $D$ are simple diagrams that are distinct up to isotopy and the $\alpha_D\neq 0\in \mathfrak{D}$ . There exists a unique simple diagram $E$ in this sum,  so that $weight(E)=weight(S)+weight(S')$, and all the other simple diagrams appearing with nonzero coefficient in the sum  have strictly lower weight. Furthermore, the coefficient $\alpha_E$ is a power of $\zeta$.   \end{theorem}

Adam Sikora informs us that working with Jozef Przytycki they proved a similar result based on a filtration coming from Dehn coordinates for the simple closed curves on the surface.  The proof is based  on the following lemma.

\begin{lemma}
Let $G$ be a four valent graph with at least one vertex, embedded in a disk $D^2$. Assume that $G$  is the union of two families of properly embedded arcs $A_1\cup A_2$ and that there are three special points $p,q,r$ in $\partial D^2$,  so that 
\begin{itemize}
\item the endpoints of the $A_1$ and $A_2$ are disjoint from one another and $\{p,q,r\}$
in $\partial D^2$, and
\item if $a_1\in A_1$ and $a_2\in A_2$ then $a_1$ and $a_2$ intersect transversely, and realize their geometric intersection number relative to their boundaries, and
\item if $a,b\in A_i$ then $a \cap b=\emptyset$ , and
\item for any arc $a\in A_1\cup A_2$, the endpoints of $a$ are separated by $\{p,q,r\}$,
\end{itemize}
Then there is an embedded triangle $\Delta$ whose sides consist of an arc of $\partial D^2$ that is disjoint from $\{p,q,r\}$, an arc contained in some $a\in A_1$ that only intersects $A_2$ in a single point which is one of its endpoints, and an arc in some  $b\in A_2$ that only intersects $A_1$ in a single point which is one of its endpoints.

We call this an {\bf outermost triangle}.  \end{lemma}

\proof The graph dissects the disk into vertices, edges and faces. The alternating sum of the numbers of vertices, edges and faces is $1$ as that is the Euler characteristic of the disk.  A face $f$ has two kinds of sides, sides in $\partial D^2$ and sides in the interior of $D^2$. Let $e_{\partial}(f)$ denote the number of sides of $f$ lying in $ \partial D^2$ and $e_{i}(f)$ the number of sides of $f$ in the interior.  Similarly, let $v_{\partial}(f)$ be the number of vertices of the face that lie in $\partial D^2$, and
$v_{i}(f)$ be the number of vertices of $f$ that lie in the interior of $D^2$.   The contribution of the face $f$ to the Euler characteristic of the disk is,
\[ c(f)=1-\frac{e_i(f)}{2}-e_{\partial}(f)+\frac{v_i(f)}{4}+\frac{v_{\partial(f)}}{2}.\]
We have that $\sum_fc(f)=1$.  The faces that are contained in the interior of the disk have an even number of sides, as their edges are partitioned into arcs of $A_1$ and arcs of $A_2$. Since the arcs of $A_1$ and $A_2$ realize their geometric intersection number the interior faces have at least four sides. Hence the largest contribution of an interior face is $0$.  A face touching the boundary can have two sides, but these faces are cut off by a single component of $A_1$ or $A_2$, and contain a point of $\{p,q,r\}$ in their boundary face by the last condition. These components can be discarded and the hypotheses of the theorem hold true for the graph formed by the smaller collection of arcs. The only remaining  faces that contribute positively to the Euler characteristic of the disk are triangles with one edge on the boundary. These contribute $\frac{1}{4}$ to the Euler characteristic. There must be at least $4$ such triangles.  That means one of those triangles does not contain a point from $\{p,q,r\}$, so it is an outermost triangle. \qed

\proof Theorem \ref{weighty}   Let $S$ and $S'$ be two simple diagrams, with associated colorings $f_S,f_{S'}:C\rightarrow \mathbb{Z}_{\geq 0}$ where $C$ is the system of proper lines defining and ideal triangulation with ideal triangles $\Delta_j$. We do not need to distinguish between embedded and folded triangles for this proof, because the combinatorial lemma above is applied in the completed components of the complement of $C$. Isotope $S$ and $S'$ so that they are transverse to one another, and the lines in $C$, and realize all geometric intersection numbers $i(S,S')$, $i(S,c)$ and $i(S',c)$ for $c\in C$. Also make sure that $S\cap S'\cap C=\emptyset$. We resolve $S*S'$ one ideal triangle at a time. The four valent graph $(S\cup S')\cap \Delta_j$ for each $\Delta_j$ satisfies the hypotheses of the lemma. To start with, $A_1$ is made up of the components of $S\cap \Delta_j$ and $A_2$ is made up of the components of $S'\cap \Delta_j$.
 Therefore  we can find an outermost triangle $\Delta\subset \Delta_j$. If we resolve the crossing of $S*S'$ at the apex of the triangle there are two resolvents. One resolvent forms a bigon with the edge of the triangle, and hence any simple diagram  descendent from this resolvent has strictly lower weight than $weight(S)+weight(S')$. 
\begin{figure}[h]\includegraphics{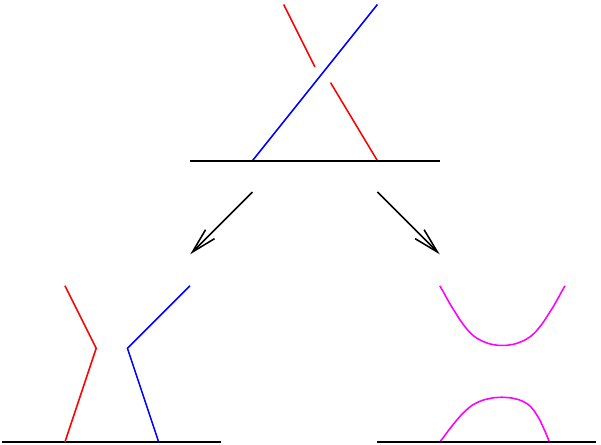}\caption{Resolving at an outermost triangle}\end{figure}

The other resolvent doesn't have a bigon.  Any state resulting in a simple diagram of weight $weight(S)+weight(S')$ is a descendent of this resolvent.  The  triangle $\Delta$ has a face $p\subset S$ and a face $q\subset S'$. Assume that $p$ lies in the component $a$ of the family $A_1$  and $q$ lies in the component $b$ of $A_2$. We smooth by forming arcs $a-p\cup q$ and $b-q\cup p$ and then perturb them slightly so that they are disjoint. To continue on inductively, we declare that the perturbed version of $a-p\cup q$ is in $A_1$, whilst removing $a$,   and the perturbate of $b-q\cup p$ is in $ A_2$ and discard $b$. Notice that the assignment  of $A_1$ and $A_2$ is now just local to the ideal triangle instead of corresponding to the diagrams $S$ and $S'$. However we work ideal triangle by ideal triangle, so this isn't a problem. If the new graph has a crossing it still satisfies the hypotheses of the lemma,
 so we can continue resolving crossings  at the apex of an outermost triangle and there is a unique resolvent  that can have a descendent of weight $weight(S)+weight(S')$. Continue until there are no crossings in $\Delta_j$. There is a  single resolvent with no bigons in $\Delta_j$ so that all the crossings in $\Delta_j$ have been resolved.  All the other resolvents with no crossings in $\Delta_j$  have bigons in $\Delta_j$ and  will lead to simple diagrams  of strictly  lower weight.  Do this for each triangle.  In the end, there is a single state of maximum weight which is a simple diagram $E$.  Since there are no bigons between $E$ and the edges of the triangulation, the admissible coloring associated to $E$ is $f_S+f_{S'}$, so $E$ is the geometric sum of $S$ and $S'$.  The  coefficient of $E$  is $\zeta^{p(E)-n(E)}$ where $p(E)$ is the number of positive smoothings and $n(E)$ is the number of negative smoothings that gave rise to the state $E$. The rest of the expansion is a linear combination of  simple diagrams  having strictly  lower weight. \qed

 If $\beta=\sum\alpha_iS_i\in K_{\mathfrak{D}}(F)$ where $\alpha_i\in \mathfrak{D}$,  and the $S_i$ are simple diagrams that are not isotopic to one another, the {\bf weight} of $\beta$ is the maximum weight of the simple diagrams $S_i$ with $\alpha_i\neq 0$.  Let
\[ \mathcal{F}_i=\{\beta\in K_{\mathfrak{D}}(F)| weight(\beta)\leq i\}.\]
It is clear that  $\mathcal{F}_i*\mathcal{F}_j\leq \mathcal{F}_{i+j}$. Filtrations like this has been used by March\'{e} \cite{M}, L\^{e} \cite{Le} and Muller \cite{Mu}.

The associated graded object $\bigoplus_i\mathcal{F}_i/\mathcal{F}_{i-1}$ is a ring extension of  $\mathfrak{D}$.  The skeins associated to $f_S:C\rightarrow \pZ$ of weight $i$ form a basis for $\mathcal{G}_i=\mathcal{F}_i/\mathcal{F}_{i-1}$. We denote the element of $\mathcal{G}_i$ associated to $S$ by $[f_S]$ where $f_S:C\rightarrow \mathbb{Z}_{\geq 0}$ is the admissible coloring corresponding to $S$.

If $S$ has weight $i$, and $S'$ has weight $j$ then the simple diagram of $f_S+f_{S'}$ has weight $i+j$.  Let $p(S,S')$ be the number of positive smoothings, and $n(S,S')$ be the number of negative smoothings associated to the state of $S*S'$ with weight $i+j$.  We define,
\[ e(S,S')=p(S,S')-n(S,S').\]  The product in $K_{\mathfrak{D}}(F)$ descends to a bilinear map,
\[ *:\mathcal{G}_i\otimes \mathcal{G}_j\rightarrow \mathcal{G}_{i+j}.\]
By Theorem \ref{weighty}, in the graded object,
\[ [f_S]*[f_{S'}]=\zeta^{e(S,S')}[f_S+f_{S'}].\]

\begin{remark} The filtration, and hence the graded object depends on the choice of ideal triangulation.  \end{remark}

If $\beta\in K_{\mathfrak{D}}(F)$, then it has a unique expression  as a finite sum $\beta=\sum_S\alpha_S S$ where $S$ are skeins corresponding to simple diagrams, with no two isotopic to one another, and $\alpha_S\in \mathfrak{D}$.   If $S_1,\ldots, S_n$ are the simple diagrams appearing with nonzero coefficient in the sum so that $weight(S_i)=weight(\beta)$, then the {\bf symbol} of $\beta$ is
\[ \sigma(\beta)=\sum_i\alpha_{S_i}[f_{S_i}]\in \mathcal{G}_{weight(\beta)}.\]

\begin{proposition} The map,
\[ \sigma:K_{\mathfrak{D}}(F)\rightarrow \bigoplus_i\mathcal{G}_i,\]
satisfies,
\[\sigma(S*(\sum_i\alpha_iS_i))=\sigma(S)*\sigma(\sum_i\alpha_iS_i)),\]
for all simple diagrams $S$ and $S_i$, so that the collection $\{S_i\}$ is linearly independent.  Furthermore,
\[weight(\beta+\gamma)\leq \max\{weight(\beta),weight(\gamma)\}.\]
If $weight(\beta)\neq weight(\gamma)$,  then
\[weight(\beta+\gamma)= \max\{weight(\beta),weight(\gamma)\}.\]
Finally , if $weight(\beta)=weight(\gamma)$, and $\sigma(\beta)+\sigma(\gamma)\neq 0$, then \[\sigma(\beta+\gamma)=\sigma(\beta)+\sigma(\gamma).\]
\end{proposition}
\qed

In \cite{Mu}, Muller defines a subtropical form that has similar properties.

\begin{lemma}  \label{nozero} Suppose that  $S_i$ is a collection  of simple diagrams with  associated admissible colorings $f_{S_i}:C\rightarrow \pZ$ all having the same weight.  Assume further that if
$f_{S_i}+f_{S_j}=f_{S_k}+f_{S_l}$ then $\{i,j\}=\{k,l\}$.
 Let $\alpha_i,\beta_i\in \mathfrak{D}$.  If, 
\[weight((\sum_i\alpha_iS_i)*(\sum_i\beta_iS_i))<weight(\sum_i\alpha_iS_i)+weight(\sum_i\beta_iS_i).\] then either all $\alpha_i=0$, or all $\beta_i=0$. \end{lemma}

\proof 

Suppose that the weight 
\[(\sum_i\alpha_iS_i)*(\sum_i\beta_iS_i)\]
is less than $weight(\sum_i\alpha_iS_i)+weight(\sum_i\beta_iS_i)$.  This means that the symbols  cancel.
Since $f_{S_i}+f_{S_j}=f_{S_k}+f_{S_l}$ if and only if $\{i,j\}=\{k,l\}$,  the cancellations in the symbol of  \[(\sum_i\alpha_iS_i)*(\sum_i\beta_iS_i)\]  can be collected as $(\zeta^{e(S_i,S_j)}\alpha_i\beta_j+\zeta^{-e(S_i,S_j)}\alpha_j\beta_i)[f_{S_i}+f_{S_j}]=0$ if $i\neq j$ and $\alpha_i\beta_i[2f_{S_i}]=0$.  If  the set $\{S_i\}$ has $n$ elements, then there are $n$ equations 
\[\alpha_i\beta_i=0\]  and
$\binom{n}{2}$ equations \[\zeta^{e(S_i,S_j)}\alpha_i\beta_j+\zeta^{-e(S_i,S_j)}\alpha_j\beta_i=0,\]
that the $\{\alpha_i\}$ and $\{\beta_i\}$ satisfy. 

The first collection of equations implies for all $i$, either $\alpha_i=0$ or $\beta_i=0$, as $\mathfrak{D}$ has no zero divisors.  Fix $i$. Without loss of generality we may assume that $\alpha_i\neq 0$. Thus $\beta_i=0$. Using this informaton in the second collection of equations,  for every  $j\neq i $, $\beta_j=0$, meaning all $\beta_j=0$. \qed

\begin{remark} A collection of skeins $\beta\in B$ spans $K_{\mathfrak{D}}(F)$ over $\mathfrak{D}$ if and only if  the set $\sigma(\beta)$ where $\beta \in B$ spans  the graded object. This is proved by   induction on the weight of a skein. \end{remark}

\begin{theorem} \label{fg} Suppose that $\mathfrak{D}$ is an integral domain and $\zeta \in \mathfrak{D}$ is a unit and $2\in \mathfrak{D}$ is a unit. Let $S_i$ be a family of simple diagrams corresponding to the integral basis of the admissible colorings of an ideal triangulation.   The skeins $\{\prod_iT_{k_i}(S_i)\}$ where the $k_i$ range over all nonnegative integers, spans $K_{\mathfrak{D}}(F)$ over $\mathfrak{D}$. \end{theorem}

\proof The symbol of $T_{k_1}(S_1)*T_{k_2}(S_2)*\ldots *T_{k_n}(S_n)$ is a power of $\zeta$ times  a simple diagram corresponding to the admissible coloring $\sum_ik_if_{S_i}$ where $f_{S_i}$ is the admissible coloring corresponding to $S_i$. Since the symbols of these skeins correspond to all simple diagrams we can inductively rewrite any skein as a linear combination of these by starting at the terms of highest weight. \qed

This extends a theorem of Bullock \cite{B}. In that paper it is proved that the arbitrary products of a finite collection of curves $S_i$ spans.  Our theorem is stronger because we can specify the order of the product of the $S_i$, as no matter what order we work in, the leading terms are the same, though maybe with different powers of $\zeta$ as the lead coefficient.  It could be that the integral basis of the space of admissible colorings is not linearly independent over $\mathbb{Z}$, so we don't have that the products form a basis.

\subsection{The case when $\zeta$ is a primitive $2N$th root of unity.}

Now we go on to study $K_N(F)$, meaning the coefficients are $\mathbb{Z}[\frac{1}{2},\zeta]$, where  $\zeta$ is a primitive $2N$th root of unity for some odd $N\geq 3$, and $A$ is set equal to  $\zeta$.  Recall, $\chi(F)$ is the image of the threading map \[\tau:K_1(F)\rightarrow K_N(F).\] The map $\tau$ threads every component of a framed link corresponding to a simple diagram with $T_N(x)$. Since \[T_k(x)=\sum_{i=0}^{\floor{k/2}}(-1)^i\frac{k}{k-i}\binom{k-i}{i} x^{k-2i},\] the symbol of $\tau(S)$ of the simple diagram corresponding to $f_S:C\rightarrow \pZ$  of weight $i$ is $[Nf_S]\in \mathcal{G}_{Ni}$.

Let $S$ be a simple diagram with associated coloring $f_S:C\rightarrow \mathbb{Z}_{\geq 0}$. Assume that $f_S$ is not identically zero. The integers $\{f_S(c)\}_{c\in C}$ generate a subgroup of $\mathbb{Z}$, which being cyclic has a smallest positive generator,
denoted $gcd(f_S)$.  

\begin{lemma}\label{divide}If $n>0$ is odd and $n|gcd(f_S)$ then $\frac{f_S}{n}:C\rightarrow \pZ$ is an admissible coloring with associated simple diagram $S'$ and $(S')^n=S\in K_N(F)$. \end{lemma}

\proof Since $F$ is orientable, the diagram $S'$ is two sided so that we can push it completely off of itself to take the product. This means that the admissible coloring  of $(S')^n$ is $nf_{S'}:C\rightarrow \pZ$.  If $f_S:C\rightarrow \pZ$ is an admissible coloring, and for all $c\in C$, the odd integer $n|f(c)$, then for any $\{a,b,c\}=\partial \Delta$ of an embedded deal triangle in the triangulation, $\{f_S(a)/n,f_S(b)/n,f_S(c)/n\}$ satisfy all three triangle inequalities as the triangle inequality is linear.  The sum $\frac{f_S(a)+f_S(b)+f_S(c)}{n}$ is even, as an even number divided by an odd number is even.  Similarly, $\frac{f_S}{n}:C\rightarrow \pZ$  satisfies the conditions to be admissible for folded triangles. \qed

\begin{proposition} Suppose that $s=\sum_j\alpha_j[f_{S_j}]\in \mathcal{G}_i$ where the $\alpha_j\in \mathbb{Z}[\frac{1}{2},\zeta]$ and the $f_{S_i}:C\rightarrow \pZ$ are distinct admissible colorings of weight $i$. The symbol $s$ is the symbol of an element of $\chi(F)$ if and only if for all $j$, $N|gcd(f_{S_j})$. \end{proposition}

\proof  If $f_{S_i}:C\rightarrow \pZ$ is the admissible coloring corresponding to $S_i$, by Lemma \ref{divide}, so is $f_{S_i}/N:C\rightarrow \pZ$.  Denote the corresponding simple diagram by $(S/N)_i$,
\[\sigma(\tau(\sum_i\alpha_i(S/N)_i))=\sigma(\sum_i\alpha_iT_N((S/N)_i))=N\sigma(\sum_i\alpha_i(S/N)_i)=\sum_i\alpha_i[f_{S_i}].\] 

On the other hand, if $S$ is a simple diagram associated to the admissible coloring \[f_S:C\rightarrow \pZ\] then the coloring associated with $\tau(S)$ is $Nf_S$, and $N$ divides $gcd(Nf_S)$.
\qed

\begin{theorem} Let $F$ be a finite type surface. If $S_i$ is any system of simple diagrams corresponding to an integral basis of the admissible colors, then the skeins $\prod_iT_{k_i}(S_i)$, where $k_i\in \{0,1,\ldots, N-1\}$ span $K_N(F)$ over $\chi(F)$. In specific $K_N(F)$ is a finite ring extension of $\chi(F)$. \end{theorem}

\proof The proof is by induction on the weight of the skein.  Start with a skein written in terms of the basis over $\mathbb{Z}[\frac{1}{2},\zeta]$ of simple diagrams,
\[ \sum_j\alpha_j \prod_iT_{k_{i,j}}(S_{i,j}),\] with $\alpha_j\in \mathbb{Z}[\frac{1}{2},\zeta]$.
Start with a term $j$ of highest weight. Since
\[T_{N+k}(x)=T_N(x)*T_k(x)-T_{|N-k|}(x)\] if some $k_{i,j}\geq N$ then  as $\chi(F)$ is central, we can
factor out an element of $\chi(F)$ from the term to get a simple diagram of lower weight. Continue on inductively till the skein is written as,
\[ \sum_j\beta_j\prod_iT_{k_{i,j}}(S_{i,j})\] where all $k_{i,j}\in  \{0,1,\ldots, N-1\}$, and $\beta_j\in \chi(F)$. 
\qed

\begin{remark} The theorem is true for closed surfaces.  If $F$ is closed and $p\in F$,
the inclusions $K_N(F-\{p\})\rightarrow K_N(F)$, and $\chi(F-\{p\})\rightarrow \chi(F)$ are surjective homomorphisms that fit into a commutative square,
\[ \begin{CD}  K_N(F-\{p\})@>>>K_N(F) \\ @AAA @AAA \\ \chi(F-\{p\}) @>>> \chi(F).\end{CD}\]
After choosing an ideal triangulation for $F-\{p\}$, if the admissible colorings associated with $S_i$ form an integral basis then $T_{k_1}(S_1)*\ldots * T_{k_n}(S_n)$ where $k_i\in \{0,\ldots, N-1\}$ span $K_N(F)$ over $\chi(F)$.

In \cite{AF} we prove that $K_N(\Sigma_{1,0})$ is not free over $\chi(\Sigma_{1,0})$, so
there are definitiely linear dependencies between the elements of the spanning set produced this way.
\end{remark}

\begin{theorem} For every $\phi:\chi(F)\rightarrow \mathbb{C}$, $K_N(F)_{\phi}$ is a finite dimensional algebra over the complex numbers. \end{theorem}  
\qed

\section{Computing the trace} 

We begin with a topological lemma.

Suppose that $C$ is a properly embedded system of disjoint lines in the finite type surface $F$.  A {\bf monogon} is a component of the complement of $C$ that completes to a closed disk with a single point removed from its boundary.

\begin{figure}[h] \includegraphics{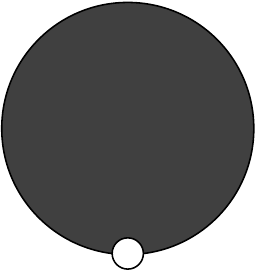}\caption{A monogon}\end{figure}

A {\bf bigon} is a component of the complement of $C$ that completes to a closed disk with two points removed from its boundary.

\begin{figure}[h] \includegraphics{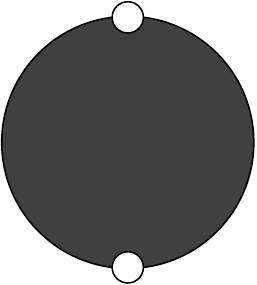}\caption{A bigon}\end{figure}

\begin{proposition}  Suppose that $C$ is a properly embedded system of disjoint lines in the finite type surface $F$ whose complement has no monogons or bigons.  There exists $D$ a collection of properly embedded lines so that $C\cup D$ defines an ideal triangulation of $F$. \end{proposition}

\qed

Suppose that $P\subset F$ is a primitive diagram. Place a complete Riemannian metric on $F$. We say that $P$ {\bf fills} $F$ if the completion of every component of $F-P$ is either a once punctured disk,  or a pair of pants.  The maps from the completions into the surface can be $2-1$ on boundary components, so the closures of the complements are possibly, a punctured disk, a pair of pants $\Sigma_{0,3}$  or a surface of genus one with one boundary component, $\Sigma_{1,1}$. 

The goal of this section is to prove that given a filling diagram $P$ which is a union of components $\{S_k\}$, there is a collection of  disjoint properly embedded  lines $\{c_j\}$ and $\alpha_{j,m}\in \mathbb{Z}_N$ so that $I_m:\{S_k\}\rightarrow \mathbb{Z}_N$ given
by $I^m(S_k)=\sum_{j}\alpha_{j,m}i(c_j,S_k)=\delta_k^m$, where $\delta_k^m$ is Kronecker's delta.  The functions  $\{I^m\}$ are ``dual'' to the components of the diagram.  We are not requiring $N$ to be prime, so $\mathbb{Z}_N$ has zero divisors, so care needs to be taken.  For the construction to work, it is essential that $N$ be odd, as there are linear dependencies between the components of $P$ as elements of $H_1(F;\mathbb{Z}_2)$, which implies that the sums of the geometric intersection numbers with components that are linearly dependent must be even. Luckily $2$ is a unit in $\mathbb{Z}_N$. The reciprocal of $2$ is $\frac{N+1}{2}$, which will be used repeatedly in the construction.

If $P$ fills $F$,  there is a dual $1$-dimensional $CW$-complex, with a $0$-cell for every component of the complement of $P$ and a $1$-cell for every component of $P$.  The trivalent $0$-cells of the $CW$-complex correspond to components of the complement that complete to pants. The monovalent $0$-cells correspond to components of the complement that complete to a punctured disk.  If a $1$-cell has both its endpoints at the same $0$-cell,  the corresponding simple closed curve is a nonseparating curve lying in the closure of a   component of the complement of $P$ that is homeomorphic to $\Sigma_{1,1}$.  The $CW$-complex minus its valence one vertices can be properly embedded in the surface $F$, where each edge intersects the corresponding simple closed curve once in a transverse point of intersection and the trivalent vertices embedded in the corresponding components of the complement of $P$, and the ends of the deleted $CW$-complex mapped to the ends of the corresponding disk with a point deleted.  The edges of the $CW$-complex are in one to one correspondence with the components of $P$, if the edge $e$ and the component $S$ intersect one another we say that they are {\bf dual}. The intersection is neccessarily a single point of transverse intersection.

 Below is a twice punctured  surface of genus three. The filling diagram is in blue, and the embedded dual graph is red.

\begin{figure}[h]\includegraphics{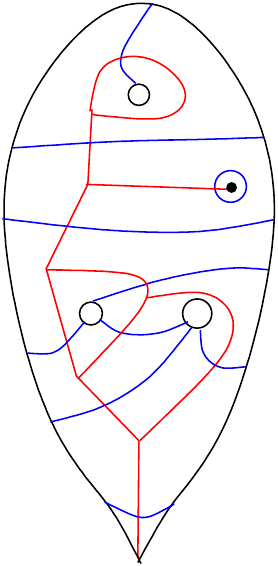}\caption{ A filling diagram and its dual graph} \end{figure}

Choose a maximal tree of the $CW$-complex  and a valence one $0$-cell. Orient the tree so that it is rooted at the chosen $0$-cell. That is, every edge is oriented so that it points towards the root.  The monovalent $0$-cells of the tree that are sources are the {\bf leaves} of the tree. The rooted tree is in red.

We will build a  train track from this tree. \begin{figure}[h] \label{tree} \includegraphics{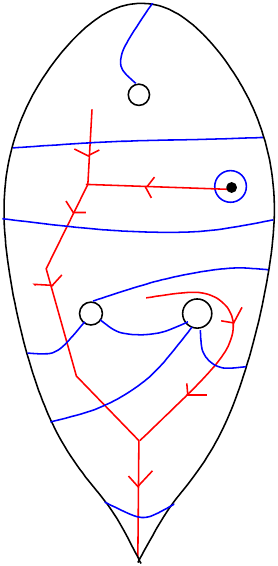}\caption{A maximal rooted tree} \end{figure}

The diagram below is color coded so that each of the following steps is visible.
First smooth the vertices ot the tree so that the two edges pointing into each interior $0$-cell have the same outward pointing tangent vector.  Next for each  component of the diagram that doesn't bound a punctured disk, that is dual to an edge of the tree, push it off itself towards the root, and then put a kink in it where it intersects the edge dual to it and smooth the kink to get a switch where both outward normals of the curve at the kink point towards the root. These are in magenta.  Next, add the remaining edges of the $CW$-complex, so that their outward normals, at the switches created, point towards the root. These are in green.  If both endpoints of the edge are attached at the same $0$-cell, that edge $e$ lies in the closure of a component of the complement of $P$ that is a torus. If $S$ is the dual edge, push it off of itself and add a kink where it intersects $e$ so that the outward tangent vectors point towards the vertex in the torus component. This is in brown.
Suppose now that the $0$-cells of the tree that $e$ is attached at are distinct. For each one of those $0$-cells that is a leaf, add a branch to the track, which is a pushoff of the dual component of $P$, with a kink in it that makes a switch in the train track pointing at that $0$-cell. These are in yellow. 

\begin{figure}[h] \includegraphics{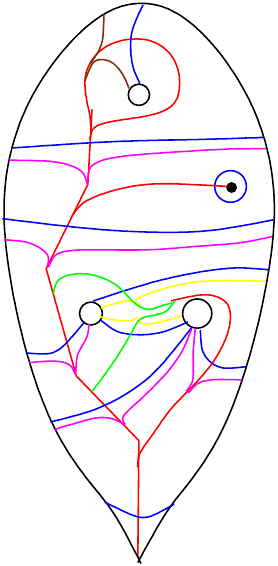} \caption{The train track}\end{figure}

We produce a family of disjoint properly embedded lines by splitting the tree at the switches and cutting open all the way to the root. The switches in the tree point towards the roots, and the switches in the additional edges point towards the tree, so the process of cutting open along switches terminates at the root, and  we have produced a family of disjoint properly embedded lines.  The train track does not carry any simple closed curves.

\begin{lemma}\label{ideal} Let $F$ be a noncompact finite surface of negative Euler characteristic, and let $C$ be a collection of lines built as above. 
There is a subset of these lines
so that for any $S\in P$ there is a linear combination \[I^S=\sum_i\alpha^i_Si(c_i, \ )\]
where the $\alpha_i\in \mathbb{Z}_N$ so that for any components $S,T$ of $P$,
$I^S(T)\equiv\delta_S^T\mod{N}$, where $\delta_S^T$ is Kronecker's delta.  This family can be seen to have no bigons or monogons in it's complement, so it can be built up to an ideal triangulation of $F$. \end{lemma}

\proof  Order the components of $P$ so that $S,T$ are dual to edges in the tree then their relative order is consistent with their distance from the root of the tree, and if they aren't dual to edges of the tree then they come after all the components that are dual to edges of the tree.  Working in order we prove that given $S$ a component of $P$ there is a line $c_S$ in our family and $\beta_S\in \mathbb{Z}_N$ so that $\beta_Si(c_S,S)\equiv 1\mod{N}$, and if $T>S$ then $c_S\cap T=\emptyset$, or we exchange order so that we can. The family $I^S$ is then produced  by taking the appropriate linear combinations of the $\beta_Si(c_S, \ )$.  Since the lines $c_S$ are indexed by $S$, the condition on intersections implies that no line is homotopically trivial (bounds a monogon) and no two lines are parallel (cobound a bigon), so the family $c_S$ can be built up to a triangulation. The complication of the construction is that to construct the line for a given edge in the tree we need to understand what immediately follows the edge in the ordering.

We start at the root. If an edge leaving the root is leaf in the tree, there are three possible cases. The surface could be a once punctured torus, or a thrice punctured sphere, or  the terminal points of the edge are at punctures, and the punctured disks containing those punctures abut the same pair of pants.  The construction for the punctured torus, and thrice punctured pair of pants can be done by inspection. We focus on the last case, shown in the figure below.

\begin{figure}\label{terminal}\scalebox{.66}{ \includegraphics{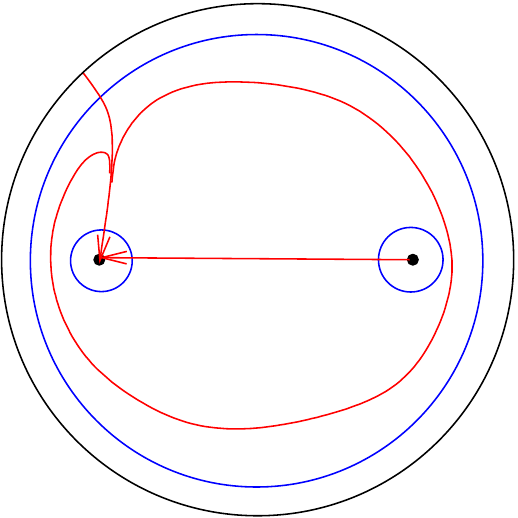}}\caption{ An edge that is simultaneously initial and a leaf}\end{figure}

According to our rules either the edge of the tree dual to the blue curve parallel to the outer boundary, or the line joining to the two punctures could come first.  You really want the edge dual to the curve parallel to the outer boundary component to be first. If $S_1$ is the component of $P$ that bounds the punctured disk at the root, let $c_{S_1}$ be the line built from the branch of the track that follows the outer boundary component before heading to the puncure. Notice $i(c_{S_1} ,S_1)=2$, so $\beta_{S_1}=\frac{N+1}{2}$ is the correct coefficient so that $\beta_{S_1}i(c_{S_1},S_1)\equiv 1 \mod{N}$.
The circle $S_2$ surrounding the other puncture has geometric intersection number one with the line $c_{S_2}$ having one end at each puncture.  Since the line $c_{S_1}$ is completely inside the diagram we have that it has geometric intersection number $0$ with all later curves.

Now suppose that the edge leaving the root is not a leaf.  In the figure below we show
the situation.  The line $c_{S_1}$ coming from the branch of the track that runs around the outer boundary component has geometric intersection number $2$ with $S_1$ and misses all the other components of the filling diagram, so letting $\beta_{S_1}=\frac{N+1}{2}$, we have that $\beta_{S_1}i(c_{S_1},\ )$ has $\beta_{S_1}i(c_{S_1},S_1)\equiv 1 \mod{N}$ and for all later components it evaluates to $0$. 

\begin{figure}[h]\scalebox{.66}{ \includegraphics{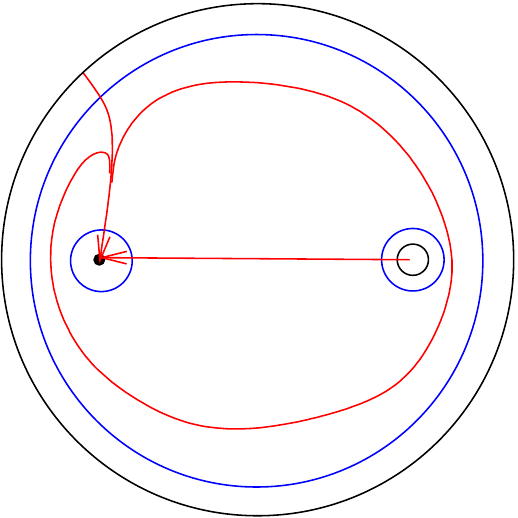}}\end{figure}

An edge dual to $S_i$ of the tree is intermediate if there is an edge dual to $S_{i-1}$ before it, and an edge dual to $S_{i+1}$ after it in the tree from the ordering. Let $c_{S_i}$ be the line coming from the branch of the track that was built by perturbing $S_{i+1}$.  Since,
$i(c_{S_i},S_i)=2$, we let $\beta_S=\frac{N+1}{2}$. Do all intermediate edges before doing the leaves. 

If an edge is a leaf, then it could end at a puncture, it could have both ends of an edge not part of the tree attached at its terminal $0$-cell, or it could have two different edges not in the tree attached at its terminal end. These both occur in Figure \ref{tree} above. The highest leaf in the diagram is of the first type, and the lower leaf is of the second kind.

In the first case, the component $S$ of $P$ bounds a punctured disk.  The line $c_S$ of the train track that emanates from that puncture and ends at the root, has geometric intersection number $1$ with $S$, we let $\beta_S=1$.

In the second case, the vertex of the edge lies in a torus that is the closure of a component of $F-P$. Call the  curve dual to the edge with both its ends attached at that $0$-cell $S'$.
The pushoff of $S'$ gives rise to the line $c_S$ that has geometric intersection $2$ number with $S$ that is dual to the edge, so $\beta_S=\frac{N+1}{2}$.

In the third case, let $S'$ be the component of $P$ that is dual to one of the edges attached at the leaf.  The pushoff of $S'$ towards the leaf gives rise to an embedded line that has geometric intersection number $2$ with the curve $S$ dual to the edge. Once again $\beta_S$ is $\frac{N+1}{2}$.   

Throw out any curves that weren't used.  Augment to form an ideal triangulation. \qed

\begin{lemma} Let $P$ be a primitive diagram with components $\{S_i\}$.  Choose an ideal triangulation as in Lemma \ref{ideal}  and let $I^j=\sum_{i}\alpha_i^ji(c_{i}, \ )$ be the linear combination of lines so that $I^j(S_k)=\delta_j^k$. Suppose that $\prod_{j\neq m}S_j^{k_j}$ is a monomial not involving $m$, then \[I^m(\prod_{j\neq m}S_j^{k_j}S_m^{k_m})\equiv k_m \mod{N}.\]

If $a \in \chi(F)$ is nonzero,  then $I^m(\sigma(a\prod_{j\neq m}S_j^{k_j}S_m^{k_m}))\equiv k_m\mod{N}$.
\end{lemma}

\proof  Since the curves $S_i$ are disjoint from one another \[I^m(\prod_{j\neq m}S_j^{k_j}S_m^{k_m})=\sum_{j\neq m}I^m(S_j^{k_j})+I^m(S_m^{k_m})=\sum_{j\neq m}k_jI^m(S_j)+k_mI^m(S_m)\equiv\]\[ \sum_{j\neq m}k_j0+k_m1\equiv k_m\mod{N}.\] 

For the second part, since $\prod_{j\neq m}S_j^{k_j}S_m^{k_m}$ is a simple diagram,
$\sigma(a\prod_{j\neq i}S_j^{k_j}S_m^{k_m})=\sigma(a)\sigma(\prod_{j\neq m}S_j^{k_j}S_m^{k_m})$.  Since $a\neq 0$ its symbol corresponds to an admissible coloring that is divisible by $N$, $g:C\rightarrow N\pZ$.  If $f_{S_j}:C\rightarrow \pZ$ is the coloring corresponding to $S_j$, the symbol of  the product $\sigma(a)\sigma(\prod_{j\neq m}S_j^{k_j}S_m^{k_m})$ is,
\[ \sum_{j}k_jf_{S_j}+g.\]  Since $i(c_{m},\sum_{j}k_jf_{S_j}+g)= \sum_jk_j(f_{S_j}(c_{m})+g(c_{m}))$, \[I^m(\sigma(a\prod_{j\neq i}S_j^{k_j}S_m^{k_m}))\equiv I^m(\prod_{j\neq m}S_j^{k_j}S_m^{k_m}))\equiv k_i \mod{N}.\]

\qed

\begin{proposition} Let $P$ be a primitive diagram with $n$ components $\{S_j\}$. Choose $S_m$ and form $\chi(F)[S_1,S_2,\ldots, \hat{S}_i,\ldots, S_n]$ meaning  adjoin all $S_j$ with $j\neq m$.  If $q(x)$ is a nonzero polynomial with coefficients in $\chi(F)[S_1,S_2,\ldots, \hat{S}_m,\ldots, S_n]$ and $q(S_m)=0$ then the degree of $q$ is at least $N$. \end{proposition}

\proof  Suppose that $p(x)=\sum_l a_lx^l$ where the $a_l\in \chi(F)[S_1,S_2,\ldots \hat{S}_m,\ldots, S_n]$.  Suppose further that $p(S_m)=0$. Finally assume that the degree of $p$ is less than $N$.  Choose an ideal triangulation as in Lemma \ref{ideal}, with the $I^k$ so that $I^k(S_j)=\delta_j^k$.  Let $\sigma(a_lS_m^l)$ be the symbol of the $l$th term of $p(S_m)$.  Since $p(S_m)=0$ the symbols cancel. However since $I^m(\sigma(a_lS_m^l))\equiv l\mod{N}$ and all the $0\leq l<N$, no two are isotopic. 
Meaning they are all zero. If the symbol is zero then so is the polynomial. \qed

If $S\in F$ is a simple closed curve, there is a polynomial of degree $N$ with coefficients in $\chi(F)$ that $S$ satisfies.  Specifically,
\[ p(x)=\sum_{i=0}^{\floor{N/2}}(-1)^i\frac{N}{N-i}\binom{N-i}{i} x^{N-2i}-T_N(S).\]

Before leaping into the final chain of theorems, we take pause
to discuss the algebraic themes involved.  They are well known, but
 bear repeating, to give shape to the argument.

\begin{proposition}\label{field} Suppose that $K$ is a field, and $p_m(x)$ with $m\in \{1,\ldots n\}$ are polynomials of positive degree with coefficients in $K$.  Suppose further that the least degree nonzero polynomial satisfied by any $S_m$ with coefficients in $K[S_1,\ldots,\hat{S_m},\ldots S_n]$ is equal to the degree of $p_m(x)$, and that $R$ is an commutative algebra over $K$ that
is generated by $S_m$, with $p_m(S_m)=0$, then if  $R$ has no zero divisors it is a  finite dimensional field extension of $K$, whose degree is the product of the degrees of the $p_m(x)$.

\end{proposition}

\proof The proof is by induction on $n$ First  the case $n=1$, the evaluation
map $K[x]\rightarrow R$ has as it's kernel a principal ideal $(q(x))$ where $q(x)|p_1(x)$.
The polynomials over a field have Krull dimension $1$ which implies prime ideals are maximal.  If $R$ has no zero divisors, $(q(x))$ is prime, so its maximal, and $R$ is a field,
and since $p_1(x)$ has minimal degree for all polynomials with coefficients in $K$ that $S_1$ satisfies $(q(x))=(p_1(x))$ and the degree of the extension is the degree of $p_1(x)$.

Assume the statement is true in the case of $n-1$ generators, and $R$ has no zero divisors, and is generated by $S_1,\ldots, S_n$. The ring $K[S_1,\ldots,S_{n-1}]$ is a subring of $R$ and has no zero divisors, so by the inductive hypothesis it is a field. If $R$ is the result of adjoining $S_n$ to the field  $K[S_1,\ldots,S_{n-1}]$, and having no zero divisors, is a field extension of $K[S_1,\ldots,S_{n-1}]$ of degree equal to $p_n(x)$. Therefore $K[S_1,\ldots,S_n]$ is a finite dimensional field extension of $K$, whose degree over $K$ is the product of the degrees of the $p_i(x)$.

\qed

\begin{proposition}\label{tensor} With the same hypotheses as Proposition \ref{field}, $R$ is naturally isomorphic to \[ \otimes_K K[S_i].\] \end{proposition}

\proof Recall that $\otimes_KK[S_i]$ is an algebra over $K$ of dimension equal to the product of the degrees of the $p_i(x)$.  The map,
\[ \psi:\otimes_KK[S_i]\rightarrow R\] given by
$\psi(S_1^{i_1}\otimes \ldots \otimes S_n^{i_n})=S_1^{i_1}\ldots S_n^{i_n}$ is an onto algebra homomorphism, which is injective because the domain and range have the same dimension. \qed

The next  proposition gives the method by which we will be computing the trace.
\begin{proposition} \label{bigtomato} Suppose that $K\leq P$ is a finite dimensional field extension and $J$ is a finite dimensional  algebra over $P$. Thus $J$ is a finite dimensional algebra over $K$.  If $s\in P$, then
it defines a $K$-linear maps $l_s:P\rightarrow P$, and $L_s:J\rightarrow J$ by left multiplication. If $d$ is the dimension of $J$ over $P$, then
\[ tr(L_s)=dtr(l_s),\] where the traces are both taken as linear maps over $K$. \end{proposition}

\proof  Since $K\leq P$ is finite dimensional it has basis $p_1,\ldots p_n$ over $K$. Since
$J$ is a finite dimensional vector space of dimension $d$ over $P$ it has basis $j_1,\ldots j_d$ over $P$.  This implies that $p_aj_c$ is a basis of $J$ over $K$.  Expressing $l_s$ respect to the basis $p_a$, we get
\[ l_s(p_a)=\sum_b l^a_b p_b.\]  Since $s$ acts as scalar multiplicationon $J$, 
\[ L_s(p_aj_c)=\sum_b l^a_bp_bj_c.\]
Hence the matrix for $L_s$ decomposes into $d$ blocks that are all copies of the matrix  for $l_s$. Therefore the trace of $L_s$ is equal to $d$ times the trace of $l_s$.\qed

\noindent Lets get to work.

\begin{proposition} \label{tech}Suppose that $S_i$ is a system of disjoint simple closed curves
on the surface of finite type, $F$,  no two of which are parallel.  That is, $\prod_iS_i$ is a primitive diagram. Then $S^{-1}\chi(F)[S_i]$ has no zero divisors.  \end{proposition}

\proof 

Choose an ideal triangulation cut out by $C$, so that there are $I^i=\sum_j\alpha^i_jc_j\in C$ with
$I^i(S_k)=\sum_j\alpha^i_j i(c_j,S_k)\equiv \delta_i^k \mod{N}$, as given by Lemma \ref{ideal}, suppose that 
\[(\prod_i\alpha_i S_i)(\prod_i\beta_iS_i)=0,\]
where $\alpha_i,\beta_i\in \chi(F)$.  (We can always clear fractions.)  Since the symbols of elements of $\chi(F)$ are divisible by $N$, if $f_{S_i}:C\rightarrow \pZ$ are the admissible colorings associated to $S_i$, there are admissible colorings $g_i,h_i:C\rightarrow \pZ$ so that
\[\sigma(\alpha_iS_i)=z_i[f_{S_i}+Ng_i]\] and \[\sigma(\beta_iS_i)=    w_i [f_{S_i}+Nh_i],\]
where $z_i,w_i\in \mathbb{Z}[\frac{1}{2},\zeta]$.

Supposing that $f_{S_i}+Ng +f_{S_j}+Nh=f_{S_k}+Nm+f_{S_l}+Np$ for any functions $g,h,m,p:C\rightarrow \pZ$.  This means that evaluating  both sides on $I^i,I^j,I^k,I^l$ will give the same answer modulo $N$. The only way this can happen is if $\{i,j\}=\{k,l\}$, so the hypotheses of Lemma \ref{nozero} are satisfied.
Therefore,
\[(\prod_i\alpha_i S_i)(\prod_i\beta_iS_i),\]
can only be zero if either the symbols of the $\alpha_i$ are all zero, which means the $\alpha_i$ are all zero or if the symbols of the $\beta_i$ are all zero which means the $\beta_i$ are all zero.  \qed

\begin{theorem} \label{dimens} If $S_i$, with $i\in \{1,\ldots,n\}$  is a system of simple closed curves on $F$ that forms a primitive diagram then $\chi(F)[S_1,\ldots,S_n]$ is a field of dimension $N^n$,
and 
\[ \chi(F)[S_1,\ldots,S_n]\cong \otimes_{\chi(F)} \chi(F)[S_i].\] \end{theorem}

\proof By Proposition \ref{tech} , we can apply Propositions \ref{field} and \ref{tensor}.\qed

\begin{theorem} \label{vs} Given  $S_i$  a system of disjoint simple closed curves
on the surface of finite type, $F$, no two of which is parallel, $S^{-1}K_N(F)$ is a finite dimensional algebra over $S^{-1}\chi(F)[S_i]$.\end{theorem}

\proof By the theorem in the last section $K_N(F)$ is a finitely generated module over $\chi(F)$. Localizing this means $S^{-1}K_N(F)$ is a finite dimensional vector space over $S^{-1}\chi(F)$.  Since $S^{-1}\chi(F)\leq S^{-1}\chi(F)[S_i]\leq S^{-1}K_N(F)$, we have that $S^{-1}K_N(F)$ is finite dimensional over $S^{-1}\chi(F)[S_i]$.\qed

If $S\subset F$ is a nontrivial simple closed curve, let $\Sigma_{0,2}(S)$ be an annular neighborhood of $S$ in $F$.  There is a left action of $K_N(\Sigma_{0,2}(S))\otimes K_N(F)\rightarrow K_N(F)$ by gluing a copy of $\Sigma_{0,2}(S)\times [0,1]$ onto the top of $F\times[0,1]$.   Notice that
it restricts to give and action $\chi(\Sigma_{0,2}(S))$ on $\chi(F)$ making $S^{-1}\chi(\Sigma_{0,2}(S))\leq S^{-1} \chi(F)$ a field extension.

\begin{remark}  It is worth mentioning that \[ S^{-1}\chi(\Sigma_{0,2}(S))[S]=S^{-1}K_N(\Sigma_{0,2}(S)).\] \end{remark}

\begin{theorem}  $S^{-1}\chi(F)[S]$ is the result of extending the coefficients of $S^{-1}\chi(\Sigma_{0,2}(S))[S]$ as a vector space over $S^{-1}\chi(\Sigma_{0,2}(S))$ to a vector space over $S^{-1}\chi(F)$. \end{theorem}

\proof  The dimension of $S^{-1}\chi(\Sigma_{0,2}(S))[S]$ over $S^{-1}\chi(\Sigma_{0,2}(S))$ is equal
to the dimension of $S^{-1}\chi(F)[S]$ over $S^{-1}\chi(F)$, so the map,
\[ S^{-1}\chi(\Sigma_{0,2}(S))[S]\otimes_{S^{-1}\chi(\Sigma_{0,2})}S^{-1}\chi(F)\rightarrow \chi(F)[S]\]
that sends $S\otimes 1$ to $S$ is a linear isomorphism. \qed

\noindent From our last paper:

\begin{proposition} \cite{AF} If $\Sigma_{0,2}$ is an annulus and $x$ is the skein at its core, and \[tr:K_N(\Sigma_{0,2})\rightarrow \chi(\Sigma_{0,2})\] is the unormalized trace,
$tr(L_{T_k(x)})=0$ unless $k|N$ at which point $tr(L_{T_{aN}(x)})=NT_{aN}(x)$.\end{proposition}
\qed

This implies the same result for $T_k(S):S^{-1}K_N(\Sigma_{0,2}(S))\rightarrow S^{-1}K_N(\Sigma_{0,2}(S))$. 

\begin{proposition} Let $S\subset F$ be a nontrivial simple closed curve. Define
$L_{T_k(S)}:S^{-1}\chi(F)[S]\rightarrow S^{-1}\chi(F)[S]$ by left multiplication, then $tr(L_{T_k(S)})=0$ unless $k|N$ at which point $tr(L_{T_{k(S)}})=NT_{k}(S)$. \end{proposition}

\proof The map $L_{T_k(S)}:S^{-1}\chi(F)[S]\rightarrow S^{-1}\chi(F)[S]$ comes from \[L_{T_k(S)}:S^{-1}\chi(\Sigma_{0,2}(S))[S]\rightarrow S^{-1}\chi(\Sigma_{0,2}(S))[S]\] by extension of scalars, and the fact that
$\chi(\Sigma_{0,2}(S))[S]=K_N(\Sigma_{0,2}(S))$. \qed

\begin{proposition}  Let $\prod_iT_{k_i}(S_i)$ act on $S^{-1}\chi(F)[S_1,\ldots,S_n]$ by multiplication \[L_{\prod_iT_{k_i}(S_i)}:S^{-1}\chi(F)[S_1,\ldots,S_n]\rightarrow S^{-1}\chi(F)[S_1,\ldots,S_n],\] the unnormalized trace of $L_{k_1,\ldots,k_n}$ is zero unless $N|k_i$ for all $i$, in which case it is $N^n\prod_iT_{k_i}(S_i)$.\end{proposition}

\proof The diagram
\[ \begin{CD}  \otimes_{\chi(F)} \chi(\Sigma_{0,2}(S_i))[S_i]\otimes_{\chi(\Sigma_{0,2}(S_i))}\chi(F) @>\psi>> \chi(F)[S_1,\ldots,S_n] \\ @V\otimes L_{T_{k_i}(S_i)}VV @VL_{\prod_iT_{k_i}(S_i)}VV  \\  \otimes_{\chi(F)} \chi(\Sigma_{0,2}(S_i))[S_i]\otimes_{\chi(\Sigma_{0,2}(S_i))}\chi(F) @>\psi>> \chi(F)[S_1,\ldots,S_n] \end{CD} \] where
$\psi$ is the natural isomorphism, commutes. This means that the trace of $L_{\prod_iT_{k_i}(S_i)}$ is the product of the traces of the \[L_{T_{k_i}(S_i)}:\chi(\Sigma_{0,2}(S_i))[S_i]\otimes_{\chi(\Sigma_{0,2}(S_i))}\chi(F)\rightarrow \chi(\Sigma_{0,2}(S_i))[S_i]\otimes_{\chi(\Sigma_{0,2}(S_i))}\chi(F)\] which are obtained by extension of scalars
from $L_{T_{k_i}(S_i)}:K_N(\Sigma_{0,2}(S_i))\rightarrow K_N(\Sigma_{0,2}(S_i))$. \qed

\begin{theorem}  Suppose that $d=[S^{-1}\chi(F)[S_1,\ldots,S_n]:S^{-1}K_N(F)]$. The unormalized trace of 
\[ L_{k_1,\ldots, k_n}:S^{-1}K_N(F)\rightarrow S^{-1}K_N(F)\] is zero unless $N|k_i$ for all $i$ in which case it is $dN^n\prod_iT_{k_i}(S_i)$. \end{theorem}

\proof By Theorem \ref{vs} $S^{-1}K_N(F)$ is a finite dimensional vector space over \[S^{-1}\chi(F)\leq S^{-1}\chi(F)[S_1,\ldots,S_n].\] so Proposition \ref{bigtomato} applies.  \qed

We define the normalized trace 
\[ Tr:S^{-1}K_N(F)\rightarrow S^{-1}\chi(F),\] to be the trace divided by $dN^n$.
The map $Tr$ is $S^{-1}\chi(F)$ linear, cyclic, and $Tr(1)=1$.

\begin{theorem}  Suppose that $s=\sum_i\beta_i P_i$ where the $\beta \in S^{-1}\chi(F)$ and the $P_i$ are primitive diagrams whose components have been threaded with $T_k$.  Let $J$ be those indices $i$ so that the components of $P_i$ have only been threaded with $T_k$ where $N|k$, then,
\[ Tr(s)=\sum_{i\in J} \beta_iP_i.\] \end{theorem}

\begin{theorem} The restriction  of $Tr:S^{-1}K_N(F)\rightarrow S^{-1}\chi(F)$ to $K_N(F)$, embedded in $S^{-1}K_N(F)$ as fractions having denominator $1$ yields,
\[ Tr:K_N(F)\rightarrow \chi(F).\]
which is a $\chi(F)$-linear map, so that $Tr(1)=1$ and for every $\alpha,\beta \in K_N(F)$,
\[ Tr(\alpha*\beta)=Tr(\beta*\alpha).\] \end{theorem}

\proof From the formula for $Tr$, the only fractions that appear in the coefficients in the trace come from fractions that are in the coefficients of the skein. \qed

\section{ The trace is nondegenerate}

\begin{lemma} Let $F$ be a finite type surface with an ideal triangulation cut out by $C$. Suppose that $\sum_iz_iS_i\in K_N(F)$, where the $z_i\in \mathbb{Z}[\frac{1}{2},\zeta]$ and the $S_i$ are distinct simple diagrams.  If for some $[f_{S}]$ appearing in the symbol of $\sum_iz_iS_i$ with nonzero coefficient $z$, has $N|gcd(f_{S})$, then \[Tr(\sum_iz_iS_i)\neq 0.\]
\end{lemma}

\proof  Suppose that the primitive diagram $P$ underlying $S$ is made up of simple closed curves $S_j'$.  The threaded diagram having lead coefficient $S$ is $\prod_jT_{Nk_j}(S'_j)$ for some $k_j\in \pZ$. Rewriting $\sum_iz_iS_i$ in terms of threaded diagrams, the threaded diagrams appearing in  the symbol appear with the same coefficients and are distinct from one another in the sum. Hence $\prod_jT_{Nk_j}(S'_j)$ appears in the trace with coefficient $z\neq 0$. This term can't cancel with other highest weight terms in the trace, as the $S_i$ were distinct nor can it cancel with lower weight terms, as that would violate the filtration of $\chi(F)$,  so $Tr(\sum_iz_iS_i)\neq 0$. \qed

\begin{theorem} Let $F$ be a noncompact, finite type surface. There are no nontrivial principal ideals in the kernel of 
\[ Tr: K_N(F) \rightarrow \chi(F).\] \end{theorem}

\proof  Suppose that $(\sum_i\alpha_i P_i)$ is a two-sided principal ideal in $K_N(F)$ where $\alpha_i\in \chi(F)$, and the $P_i$ are primitive diagrams whose components have been threaded with Tchebychev polynomials of the first type where the threadings come from $\{0,\ldots, N-1\}$. Also assume that 
$\sum_i\alpha_i P_i\neq 0$, and the diagrams $\{P_i\}$ are a linearly independent set over $\chi(F)$.   
Rewrite the sum as $\sum_jz_jS_j$ where the $z_j\in \mathbb{Z}[\frac{1}{2},\zeta]$, and the $S_j$ are distinct simple diagrams.  Choose an ideal triangulation $C$.  As the skein $\sum_jz_jS_j\neq 0$,
it's symbol is nonzero.  If some $f_S:C\rightarrow \pZ$ appearing in its symbol with nonnegative coefficient has $N|gcd(f_S)$ we are done by the last lemma. If
not, choose $f_S:C\rightarrow \pZ$ appearing in the symbol with nonzero coefficient.
We can write $S=\prod_mJ_m^{a_mN+r_m}$ where the $J_m$ are a disjoint system of simple closed curves, $a_m\in \pZ$ and $r_m\in \{0,\ldots,N-1\}$.  Consider the simple diagram,
\[ \prod_mJ_m^{N-r_m}.\]
Since $(\sum_i\alpha_i P_i)$ is an ideal, 
\[  \prod_mJ^{N-r_m}*\sum_jz_jS_j\in (\sum_i\alpha_i P_i).\]
If $\sum_iz_i[f_{S_i}]$ is the symbol of $\sum_jz_jS_j$, then the symbol of the product
is
\[ \sum_iz_i\zeta^{e(\prod_mJ^{N-r_m},S_i)}[f_{S_i}+f_{\prod_mJ^{N-r_m}}].\]
The coefficient of $[f_{S}+f_{\prod_mJ^{N-r_m}}]$ is nonzero, and $N|gcd(f_{S}+f_{\prod_mJ^{N-r_m}})$ as the corresponding diagram is,
\[ \prod_mJ_m^{(a_m+1)N}.\]  By the last lemma, the ideal $(\sum_i\alpha_iP_i)$ contains an element with nonzero trace. \qed

\begin{corollary} $S^{-1}K_N(F)$ is a symmetric Frobenius algebra over $S^{-1}\chi(F)$.\end{corollary}

\qed

\begin{corollary} There is a proper subvariety of the character variety of $\pi_1(F)$ away from which $K_N(F)_{\phi}$ is a Frobenius algebra over $\mathbb{C}$. \end{corollary}

\qed

\end{document}